\documentclass[10pt,a4paper,twoside,reqno]{amsart}

\usepackage{amsfonts}
\usepackage{amssymb}
\usepackage{amsmath}    %multi-line etc. Subsumes ams-text,amsbsy,amsopn
\usepackage{amsthm}     %provide proof env
\usepackage{amscd}      %commutative diagrams
\usepackage{euscript}
\usepackage[english]{babel}

\newcommand{\NN}{{\mathbb N}}
\newcommand{\ZZ}{{\mathbb Z}}

\newcommand{\QQ}{{\mathbb Q}}
\newcommand{\CC}{{\mathbb C}}

\newcommand{\PP}{{\mathbb P}}
\newcommand{\GG}{{\mathbb G}}

\newcommand{\PPn}{{\PP^n}}

\newcommand{\OPN}{{\mathcal O}_{\PPn}}

\DeclareMathOperator{\HH}{H}
\DeclareMathOperator{\e}{\chi}
\DeclareMathOperator{\hh}{h}
\DeclareMathOperator{\Ext}{Ext}
\DeclareMathOperator{\Hom}{Hom}
\DeclareMathOperator{\Aut}{Aut}

\DeclareMathOperator{\SL}{SL}
\DeclareMathOperator{\GL}{GL}

\DeclareMathOperator{\Id}{Id}
\DeclareMathOperator{\Stab}{Stab}
\DeclareMathOperator{\codim}{codim}
\DeclareMathOperator{\rank}{rank}

\DeclareMathOperator{\diag}{diag}
\DeclareMathOperator{\Hilb}{Hilb}

\newtheorem{thm}{Theorem}[section]
\newtheorem{lemma}[thm]{Lemma}
\newtheorem{prop}[thm]{Proposition}

\newtheorem{corol}[thm]{Corollary}
\newtheorem{remark}[thm]{Remark}
\newenvironment{proo}{{\bf Proof.}}{{\hfill $\square$}}

\title{On a compactification of the moduli space of the rational normal curves }
\author{Paolo Cascini}
\address{Dipartimento di Matematica\\ Viale Morgagni 67 A\\ 50134 Firenze\\ Italy}
\email{cascini@math.unifi.it}
\subjclass{14F05}
\keywords{moduli space, rational normal curve, vector bundle.}

\begin{document}

\begin{abstract}
For any odd $n$, we construct a smooth  minimal (i.e. obtained by
adding an irreducible hypersurface)  compactification $\mathcal
M_n$ of the quasi-projective homogeneous variety
$S_{n}=\PP\GL(n+1)/\SL(2)$ that parameterizes  the  rational
normal curves in $\PP^n$. $\mathcal M_{n}$ is isomorphic to a
component of the Maruyama scheme of the semi-stable sheaves on
$\PP^n$ of  rank $n$ and Chern polynomial $(1+t)^{n+2}$. This will
allow us to explicitly compute the Betti  numbers of $\mathcal
M_n$.

In particular $\mathcal M_{3}$ is isomorphic to the variety of
nets of quadrics defining twisted cubics, studied by G.
Ellinsgrud, R. Piene and S. Str{\o}mme \cite{EPS}.

\end{abstract}

\maketitle

\tableofcontents

\vskip 1 cm

%%%%%%%%%%%%%%%%%%%%%%%%%%%%%%%%%%%%%%%%%%%%%%%%%%%%%%
\section{Introduction}
%%%%%%%%%%%%%%%%%%%%%%%%%%%%%%%%%%%%%%%%%%%%%%%%%%%%%%

A rational normal curve $C_n$, or equivalently a Veronese curve, is a
smooth, rational, projective 
curve of degree $n$, in the complex projective space $\PP^n$: in particular
the Hilbert polynomial of $C_n$ is $P_{C_n}(d)=nd+1$. For a description of some
interesting properties of this curve, see \cite{Harris}.

The set $S_n$ of the rational normal curves is an homogeneous
quasi-projective variety isomorphic to $\PP\GL(n+1)/\SL(2)$. The
purpose of the paper is to describe a nice compactification of such
variety, by considering some vector bundles on $\PP^n$, called
Schwarzenberger bundles \cite{Sch}.
In particular we compute the Euler characteristic of such
compactification  and its Betti numbers.

\vskip 5  mm
There are several ways to define a compactification of the variety
$S_n$: probably the most natural way is to consider the
closure $\mathcal H_n$ of the open sub-scheme  of the Hilbert scheme
$\Hilb^{P_{C_n}}(\PP^n)$, parameterizing the rational normal curves in
$\PP^n$.
In \cite{PS}, the authors describe such  compacti\-fication in the case
$n=3$. In particular, they show that ${\mathcal H_3}\subseteq
\Hilb^{3d+1} (\PP^3)$ is a smooth irreducible variety
of dimension $12$. Only recently, it was proven by M. Martin-Deschamps
and R. Piene \cite{MP} that
$\mathcal H_4$ is singular. Moreover it is not difficult to verify,
with the help of the algorithm described in \cite{NS},  that $\mathcal
H_5$ and $\mathcal H_6$ are singular in the points represented by 
the $5-fold$ and $6-fold$ lines respectively. Therefore we can suspect
that $\mathcal H_n$ is singular for any $n\ge 4$ (see also \cite{Kap},
remark 2.6).

Another natural compactification is given by the closure $\mathcal
C_n$ of the quasi projective variety $S_n$ considered as an
open subset of the Chow variety $\mathcal C_{1,n}(\PP^n)$ that parameterizes
the effective cycles of dimension $1$ and degree $n$ in $\PP^n$.

In \cite{EPS}, a third natural compactification $\mathcal M_n$ of
$S_n$ is described: 
this is made by considering the space of all the $2\times n$ matrices
with linear forms as entries. In fact all the rational normal curves
in $\PP^n$ is the zero locus of the $2-$minors of such a matrix.
In particular, when $n=3$, $\mathcal M_3$ can be seen as the variety
parametrizing the 
nets of quadrics in $\PP^3$ and $\mathcal H_3$ is the blow-up of $\mathcal
M_3$.

In \cite{C}, it is shown that for any odd  $n$, the projective
variety $\mathcal M_n$ is isomorphic to a smooth irreducible
component of the Maruyama scheme  $\mathcal M_{\PP^n}
(n;c_1,\dots,c_n)$ parameterizing the semi-stable sheaves on
$\PP^n$ of rank $n$ and with Chern polynomial $c_t=\sum c_it^i =
(1+t)^{n+2}$. $\mathcal M_n$ can be seen as the  quotient of a
projective space $\PP^N$, by the action of a reductive algebraic
group $G$. This description will allow us to apply a technique of
Bialynicki-Birula \cite{B}, to compute the Betti numbers of
$\mathcal M_n$ (see also \cite{ES}).

\vskip 1 cm

Throughout the paper we will use the following notations:
\begin{itemize}

\item $V$, $W$, $I$ are complex vector spaces of dimension $n+1$,
$m+k$ and $k$ respectively, where $m\ge n$.

\item For any $A\in\PP(\Hom(W,V\otimes I))$, the cokernel $\mathcal
F_A$ of the associated map $$A^*: I\otimes\mathcal
O_{\PP(V)}\longrightarrow W\otimes\mathcal O_{\PP(V)}(1)$$ is a
coherent sheaf of rank $m$. If $A^*$ is injective and $\mathcal
F_A$ is a vector bundle, then it is said Steiner bundle of rank
$m$, and it is contained in the exact sequence:
\begin{equation}\label{suc.esatta}
0\longrightarrow I\otimes\mathcal O_{\PP(V)}\stackrel
{A^*}\longrightarrow W\otimes\mathcal O_{\PP(V)}(1)\longrightarrow
\mathcal F_A\longrightarrow 0.
\end{equation}
Moreover if $k=2$ and $n=m$, then all the Steiner bundles are
Schwarzenberger bundles (see also \cite{C}).

\item $\GG(k,n+1)$ ($\simeq\GG(k-1,\PP^n)$) is the Grassmannian of the
$k$-subspaces of $V$ or equivalently of the $k-1$ subspaces of the
projective space  $\PP^n$.

\item Let $G=\SL(I)\times\SL(W)$ and $X=\PP(\Hom(W,V\otimes I))$:
we will study the natural action of $G$ on $X$ and we will denote
by $X^s$ (resp. $X^{ss}$) the open subset of the stable (resp.
semi-stable) points of $X$.

\item For any $A\in X$, $\Stab_G(A)=\{(P,Q)\in G|PAQ^{-1}=kA \text{
    for some } k\in \CC^*\}$ is the stabilizer of $A$ by the group
  $G$.

\item $\mathcal M_{n,m,k}=X^{ss}//G$ (resp. $X^{s}/G$) is the
categorical (resp. geometric) quotient of  $X$ by $G$. In
particular, if $n=m$, we will denote $\mathcal M_{n,k}=\mathcal
M_{n,n,k}$.

\item $V^*=\CC[x_0,\dots,x_n]_1$ is the dual space of $V$.

\item For any $A\in X$, $D(A)$ is the degeneracy locus of
$A$ and $D_0(A)$ is the variety of all the points $x\in \PP^n$
such that $\rank A_x=0$.

\item ${\mathcal S}=\{ A\in X| D(A)=\emptyset\}=\{A\in X|S(\mathcal
F_A)=\emptyset \}\subseteq X^s$.

\item $S_{n,m,k} = {\mathcal S}/G$ is the moduli space
of the rank $m$ Steiner bundles on $\PP^n=\PP(V)$: in particular
$S_{n,k}=S_{n,n,k}$ is the moduli space of the ``classical''
Steiner bundles or rank $n$ on $\PP(V)$.

\item For any matrix $A\in {\mathcal M}(k\times
(m+k),V^*)$, if $A=(a_{i,j})$%_{i=0,\dots,k-1}^{j=0,\dots,m+k-1}$
we define $i_s(A)=\min\{j=0,\dots,n+k-1|a_{s,j}\neq 0\}$ (we will
often write $i_s$ instead of $i_s(A)$).

\item $j(n)=[\frac {n+3} 2]$ where $[m]$ denotes the integer part of $m$.

\item For any coherent sheaf $\mathcal E$ of rank $r$ on $\PP^n$ and for any $t\in
  \ZZ$, we write $\mathcal E(t)$ instead of $\mathcal E\otimes
  \OPN(t)$. $\mathcal E_N$ will denote the normalized of $\mathcal
  E$, i.e. $\mathcal E_N=\mathcal E(t_0)$ where $t_0\in\ZZ$ is such that
  $-r< c_1(\mathcal E(t_0))\le 0$.

\noindent Moreover, we define the slope of $\mathcal E$ as the
number $\mu(\mathcal E)=\frac {c_1(\mathcal E)} r$
  and $\mathcal E$ is said to be $\mu$-stable if it is Mumford-Takemoto stable.

\end{itemize}

\vskip 1 cm

In the first part of the paper we describe the (semi-)stable
points of the projective space $\PP(\Hom(W,V\otimes I))$ under
the action of $\SL(I)\times\SL(W)$ (see \cite{MFK} for an
introduction to the geometric invariant theory) and in particular we
will prove that, if $m<\frac{n k} {k-1}$, then   all the Steiner
bundles are defined by stable matrices,
i.e. $S_{n,m,k}\subseteq\PP(\Hom(W,V\otimes I))^s/(\SL(I)\times\SL(W))$.

\vskip 5 mm

In the second part of the paper, we investigate some properties of
$S_{n,m,2}$ and in particular of $S_{n,2}$, the moduli space of the
Schwarzenberger bundles.
By the previous correspondence of bundles and curves, $\mathcal
M_{n,n,2}$ gives us a compactification of the set of the
rational normal curves in $\PP^n$.

We define a filtration of  $\mathcal M_{n,m,2}$ and we show
that the compactification is obtained by adding an irreducible
hypersurface.

Moreover in \cite{C} it is shown that, if $k=2$ and $m$ is odd,
then $A\in\PP(\Hom(W,V\otimes I)$
 is stable if and only if the correspondent coherent sheaf $\mathcal F_A$ is
$\mu-$stable. This yields the theorem:
\begin{thm}\label{main.thm}
$M_{n,m,2}$ is isomorphic to the connected component of the Maru\-yama
moduli space  $\mathcal
M_{\PP^n}(m,c_1,\dots,c_n)$ containing the Steiner bundles. Such
component is smooth and irreducible.
\end{thm}

\vskip 4 mm

In the last two sections we compute the Betti
and Hodge numbers of the smooth projective variety $\mathcal
M_{n,m,2}$. This formula will be obtained by studying a natural action of
$\CC^*$ on $\mathcal M_{n,m,2}$: in particular we will describe its
fixed points and we will compute the weights of the action of $\CC^*$
induced in the tangent spaces of the variety at the fixed points.

\vskip 2 cm

%%%%%%%%%%%%%%%%%%%%%%%%%%%%%%%%%%%%%%%%%%%%%%%%%%%%%%
\section{The categorical quotient of $\PP(\Hom(W,V\otimes I))$ by
$\SL(I)\times\SL(W)$}
%%%%%%%%%%%%%%%%%%%%%%%%%%%%%%%%%%%%%%%%%%%%%%%%%%%%%%

We are interested in the study of the action of $G=\SL(I)\times\SL(W)$
on the projective space $X=\PP(\Hom(W,I\otimes V))$. In fact, as shown
in the introduction, each $A\in X^{ss}$, such that $A^*: I\otimes
\OPN\rightarrow W\otimes \OPN(1)$ is injective, corresponds to a
coherent sheaf $\mathcal F_A$ contained in the exact sequence
(\ref{suc.esatta}).

Furthermore $\mathcal F_A\simeq \mathcal F_B$ if and only if $P A=B Q$
for some $P\in \SL(I)$ and $Q\in \SL(W)$ (see for instance \cite{AO} or
\cite{MT}).

\

\begin{lemma}\label{iniettivo}
Let $A\in X^{ss}$.
Then both $A:W\to I\times V$ and $A^*:I\to W\times V$ are injective.
\end{lemma}
\begin{proo}
Let $A:W\to I\times V$ be non-injective. Then we can suppose that the
first column of $A$ is zero. Let us consider the 1-dimensional
parameter subgroup $\lambda:\CC^*\rightarrow G$
defined by $t\mapsto (\Id,\diag(t^{-(m+k-1)},t,\dots,t))\in
\SL(I)\times\SL(W)$:
then $\lim_{t\to 0}\lambda(t)A = 0$ and,  by the
Hilbert-Mumford criterion, the matrix $A$ cannot be semi-stable.

Let us suppose now that $A^*:I\to W\times V$ is not injective:
i.e. the first row of $A$ is zero. In this case it
suffices to consider the 1-dimensional parameter subgroup
$\mu: t\mapsto (\diag(t^{-(k-1)},t,\dots,t),\Id)\in\SL(I)\times\SL(V)$
in order to have $\lim_{t\to 0}\mu(t)A = 0$.
\end{proo}

\

As a direct consequence of the lemma, it follows that for any $A\in
X^{ss}$, the sheaf $\mathcal F_A$ is well-defined as the cokernel of
$A^*$ and is contained in the sequence (\ref{suc.esatta}). Moreover it
results $T_A:= A(W)\in \GG(m+k,I\otimes V)$.
Thus, in order to study the (semi-)stable point of $X$
by the action of $G$, it suffices to study the action of $\SL(I)$ on
the variety $\GG(m+k,I\otimes V)$: in particular we have that the
categorical quotient $\mathcal M_{n,m,k}:= X^{ss}//G$ is isomorphic to
the quotient $\GG(m+k, I\otimes V)^{ss}//\SL(I)$.

Let us recall first the following known result:
\begin{prop}
Let $T\in \GG(m+k,I\otimes V)$. The following are equivalent:
\begin{enumerate}
\item $T$ is semi-stable (resp. stable) under the action of $\SL(I)$;
\item for any non-empty subspace $I'\subsetneq I$
$$\frac {\dim T'} {\dim I'} \le \frac {\dim T} {\dim I}
\quad (\text{resp. }<)$$
where $T'=(I'\otimes V)\cap T$.
\end{enumerate}
\end{prop}
\begin{proof}
See for instance \cite{NT} (prop. 5.1.1)
\end{proof}

As a corollary we get a description of the (semi-)stable
points of $X$ by the action of $G$:

\begin{thm}\label{stabile}
$A\in X$ is not stable under the action of $G$ if and only if
with respect to suitable bases of $W$ and $I$, it results
$i_0(A)\ge i_1(A)\ge\dots\ge i_{k-1}(A)$ and there exists $s\in\{0,\dots,k-1\}$
such that:
\begin{equation}\label{condizione.stabile}
\text { either } \quad i_s(A)\ge \frac {m+k} k (k-1-s) \text{ if }
s\neq k-1 \qquad \text { or } \quad  i_{k-1}(A)>0
\end{equation}
\end{thm}

\begin{thm}\label{semistabile}
$A\in X$ is not  semi-stable under the action of $G$ if and only if
with respect to suitable bases of $W$ and $I$,
it results $i_0(A)\ge i_1(A)\ge\dots
i_{k-1}(A)$ and there exists $s\in\{0,\dots,k-1\}$
such that:
\begin{equation}\label{condizione.semistabile}
i_s(A)> \frac {m+k} k (k-1-s)
\end{equation}
\end{thm}

\vskip 1 cm

\begin{corol}\label{relat.primi}
$X^s=X^{ss}$ if and only if $(m,k)=1$
\end{corol}

\begin{proo}
If there exists $A\in X$ properly semi-stable, then there exists $s\in
\{0,\dots,k-2\}$ such that
$$i_s(A)= \frac {m+k} k (k-1-s).$$
Since $1\le k-1-s\le k-1$, such $s$ exists if and only if $(m,k)\neq 1$.
\end{proo}

\vskip 2 cm

Now we are interested to study the stability of the matrices defining
the Steiner bundles and thus we will consider all the matrices $A$
such that $\rank A_x = k$ for any
$x\in\PP^n$: in
\cite{AO} it is shown that if $n=m$ (boundary format) then all such
matrices are stable.
We generalize such result with the following:

\begin{thm}\label{vect.bundle} If $m<\frac{nk} {k-1}$ then every
indecomposable vector bundle $\mathcal F_A$ is defined by a
G.I.T. stable matrix $A$.
\end{thm}

Before proving the theorem, we remind the following known lemma:
\begin{lemma}\label{thm2.8}
Let $F$ be a vector bundle of rank $f$ on a smooth projective variety
$X$ such that $c_{f-k+1}(F)\neq 0$ and let $\phi:\mathcal
O_X^k\longrightarrow F$ be a morphism with $k\le f$. Then the
degeneracy locus $D(\phi)=\{x\in X|\rank (\phi_x)\le k-1\}$ is
nonempty and $\codim D(\phi)\le f-k+1$.
\end{lemma}

\vskip 5 mm

\begin{proof}[proof of theorem \ref{vect.bundle}]
Let $\mathcal F_A$ be an indecomposable vector bundle. Then for any
base of $W$ and $I$, $i_{k-1}(A)=0$, otherwise $\mathcal F_A=
\mathcal F'\oplus \OPN(1)$ for some vector bundle $\mathcal F'$.

Let  $I'\subseteq I$ of dimension $r$: if $s = \dim(I'\otimes V)\cap
T_A$ and $I''\subseteq I$ is such that
$I'\oplus I''=I$, then the restriction of $A^*$ in $I''$ defines a
morphism of vector bundles $A':\OPN^{k-r}\longrightarrow \OPN(1)^{m+k-s}$.

Let us suppose $s>m-n+r$, then
$c_{(m+k-s)-(k-r)+1}(\OPN(1)^{m+k-s})\neq 0$: lemma \ref{thm2.8} implies that
the degeneracy locus of $A'$ is not empty, which leads to  a
contradiction.

Thus:
$$\dim(I'\otimes V)\cap T_A\le m+k-n-k+r = m-n+r;$$
and in particular, if $m<\frac{nk} {k-1}$, it results $\dim(I'\otimes V)\cap
T_A<\frac{r(m+k)} k$, i.e. $A$ is G.I.T. stable.
\end{proof}

\begin{remark}
By lemma \ref{iniettivo} we have that if $A\in X^{ss}$, then
$A:W\hookrightarrow I\otimes V$ is injective, thus it results
$X^{ss}=\emptyset$ if $m>kn$. Furthermore it is easy to see that if
$m=nk$ the only point of $M_{n,kn,k}$ is represented by the vector
bundle $I\otimes T_{\PP^n}$.

\end{remark}

\vskip 2 cm

%%%%%%%%%%%%%%%%%%%%%%%%%%%%%%%%%%%%%%%%%%%%%%%%%%%%%%
\section{Compactification of $S_{n,m,2}$}
%%%%%%%%%%%%%%%%%%%%%%%%%%%%%%%%%%%%%%%%%%%%%%%%%%%%%%

So far we have studied the G.I.T. compactification of $S_{n,m,k}$ for
any value of $n,m$ and $k$.

From now, we restrict our study to the case $k=2$: in particular we
know that the moduli space $S_{n,2}$ is uniquely composed by Schwarzenberger
bundles and thus it is isomorphic to $$\PP\GL(n+1)/\SL(2).$$ Hence
$\mathcal M_{n,2}$ is a compactification of the set of rational normal
curves in $\PP^n$.

After a short review of the previous section, we define a
$G-$invariant filtration of the space $\mathcal M_{n,m,2}$ and we study
some properties of it.

\

Theorems
\ref{stabile} and \ref{semistabile} become:
\begin{thm}\label{stabile.k2}
Let $j(m)=[\frac {m+3} 2]$. $A\in X$ is not stable if and only if
$$\text{either}\quad A\sim\begin{pmatrix}0 &\dots& 0 & f_{j(m)+1} &\dots&
f_{m+2}\cr g_1 &\dots&  g_{j(m)} & g_{j(m)+1} &\dots& g_{m+2}
\end{pmatrix} \quad \text{or} \quad A\sim \begin{pmatrix} 0 *
  \dots * \cr 0 * \dots * \end{pmatrix}$$
\end{thm}

\vskip 4 mm

\begin{thm}\label{semistabile.k2}
If $n$ is odd then $X^{ss}=X^{s}$, i.e. there are not properly
semi-stable points in $X$.
If  $n$ is even then $A\in X$ is not semi-stable if and only if
$$\text{either}\quad A\sim\begin{pmatrix}0 &\dots& 0 & f_{j(m)+2} &\dots&
f_{m+2}\cr g_1 &\dots&  g_{j(m)+1} & g_{j(m)+2} &\dots& g_{m+2}
\end{pmatrix} \quad \text{or} \quad A\sim \begin{pmatrix} 0 *
  \dots * \cr 0 * \dots * \end{pmatrix}$$
\end{thm}

\vskip 1 cm

\begin{lemma}\label{lemma1}
Let $m$ be even and for any $i=1,2$ let us define the subspaces
  $I^i_f=<f_0^i\dots f^i_{\frac m 2}>$
  and $I_g=<g^i_0\dots g^i_{\frac m 2}>$
  of $\CC[x_0,\dots,x_n]$ of dimension
   ${\frac m 2}+1$. Moreover let
$$A^i=\begin{pmatrix} 0&\dots &0& f^i_0&\dots& f^i_{\frac m 2}\cr g^i_0&
\dots& g^i_{\frac m 2}& 0& \dots& 0  \end{pmatrix}\qquad i=1,2$$
Then
\begin{equation}\label{s1}
A^1\sim A^2
\end{equation}
if and only if
\begin{equation}\label{s2}\text{either}\quad
 I^1_f=I^2_f \text{ and } I^1_g=I^2_g \quad \text{or}\quad
 I^1_f=I^2_g \text{ and } I_g^1=I_f^2
\end{equation}
\end{lemma}

\begin{proof}
Let us suppose that (\ref{s1}) holds, then $A^1$ and $A^2$ have the same
degeneracy locus and this implies:
$$V(I^1_f)\cup V(I^1_g) = V(I^2_f)\cup V(I^2_g).$$
Since $V(I^1_f)$, $V(I^1_g)$, $V(I^2_f)$ and $V(I^2_g)$ are
irreducible, it results $V(I^1_f)=V(I^2_f)$ and $V(I^1_g)=V(I^2_g)$
or $V(I^1_f)=V(I^2_g)$ and $V(I^1_g)=V(I^2_f)$, thus (\ref{s2}) holds.

Vice-versa let us suppose $I^1_f=I^2_f$ and $I^1_g=I^2_g$ and let
$B_1,B_2\in\SL(\frac m 2 +1)$ be the respective base change matrices.
Then
$$A_1\begin{pmatrix} B_2 &
   0\cr 0 & B_1\end{pmatrix}= A_2.$$

Otherwise if $I_f=I_g'$ and $I_g=I_f'$ then if $C_1,C_2\in\SL(\frac m 2
+1)$ are the respective base change matrices, then
$$\begin{pmatrix} 0 & 1\cr 1 & 0\end{pmatrix}A_1
   \begin{pmatrix} 0 & C_2
   \cr C_1 & 0\end{pmatrix}= A_2.$$
Thus (\ref{s1}) holds.
\end{proof}

\vskip 1 cm

\begin{thm}\label{sottoinsieme}
Let $m$ be even. Then
\begin{equation}\label{sottoinsieme.formula}
(X^{ss}\setminus X^s)//G\simeq
  S^2\GG\left({\frac m 2},\PP(V)\right).
\end{equation}
\end{thm}

\begin{proof}
Let $A\in X^{ss}\setminus X^s$. Then
$$A\sim
\begin{pmatrix}0 &\dots& 0 & f_{j(m)+1} &\dots&
f_{m+2}\cr g_1 &\dots&  g_{j(m)} & g_{j(m)+1} &\dots& g_{m+2}
\end{pmatrix}$$ and thus if we consider the 1-dimensional parameter
subgroup defined by the weights $\beta=(-1,1)$ and
$\gamma=(-1,\dots,-1,1,\dots,1)$, it results:
$$\lim_{t\rightarrow 0}tA =
\begin{pmatrix}0 &\dots& 0 & f_{j(m)+1} &\dots&
f_{m+2}\cr g_1 &\dots&  g_{j(m)} & 0 &\dots& 0\end{pmatrix}.$$

Thus the points of $(X^{ss}\setminus X^s)//G$ are in one-one
correspondence with the orbits
of the matrices $\begin{pmatrix} 0&\dots&0&*&\dots&*\\
  *&\dots&*&0&\dots&0\end{pmatrix}\in X^{ss}$ by the action of $G$.
The previous lemma implies the isomorphism in
(\ref{sottoinsieme.formula}).
\end{proof}

\vskip 1 cm

\begin{remark} Since $\GG(m+2,2(n+1))\simeq \GG(2n-m,2(n+1))$, it
follows that $\mathcal M_{n,m,2}\simeq \mathcal M_{n,2n-m-2,2}$. In
particular $\mathcal M_{n,2}$ parameterizes the $n\times 2$ matrices
with entries in $V^*$: in fact a rational normal curve is the zero
locus of the minors of such a matrix.

In the case $n=3$, we have that $\mathcal M_{3,2}$ is
isomorphic to the variety of the nets of quadrics that define the
twisted cubics in
$\PP^3$. In \cite{EPS}, the authors describe this variety and they
show that there exists a natural morphism from the Hilbert scheme
compactification $\mathcal H_3$ to $\mathcal M_{3,3,2}$. It would be
interesting to know if there exist a canonical  morphism, $\mathcal H_n\to
\mathcal M_{n,n,2}$, for any odd $n$.
\end{remark}

\vskip 1 cm

For any $\omega \in I$ we define $R_\omega=\omega\otimes V\subseteq
I\otimes V$: by theorems \ref{stabile.k2} and \ref{semistabile.k2} we
have that an injective matrix $A: W\hookrightarrow I\otimes V$ is
semi-stable (resp. stable) if and only if
$$\dim R_\omega\cap T_A\le \frac{m+2} 2 \quad (\text{resp.} <)$$
for any $\omega\in I$.

For any $j=0,1,\dots$ we construct the subsets:
\begin{eqnarray*}
S^j &=&\left\{A\in X^{ss}|\exists ~\omega\in I \text{ such that } \dim
R_\omega\cap T_A\ge j+m-n \right\}\subseteq X^{ss}\quad\text{and}\\
\tilde S^j&=&\{A\in X^{ss}|\dim D(A) \ge j - 2\}\subseteq X^{ss}.
\end{eqnarray*}
Such subsets of $X$ define two filtrations:
\begin{eqnarray*}
\emptyset= S^{j_0+1}\subseteq & S^{j_0}\subseteq \dots \subseteq
&S^2\subseteq S^1=X^{ss}\\
\emptyset\subseteq\dots\subseteq \tilde S^{j_0+1}\subseteq &\tilde
 S^{j_0}\subseteq \dots\subseteq &\tilde S^2\subseteq \tilde S^1=X
\end{eqnarray*}
where $j_0=j(m)+n-m$. It results
$S^{j_0}=X^{ss}\setminus X^s$ and in particular it is empty if $m$ is odd.
Furthermore we have:
\begin{thm}\label{filtrazione}
\
\begin{enumerate}
\item $S^j\subseteq\tilde S^j\subseteq S^{j-1}$ for any $j\ge 2$;
\item $S^{2}=\tilde S^{2}$;
\item $S^{1}=\tilde S^{1}=X^{ss}$.
\end{enumerate}
In particular such subsets define a unique filtration $G-$invariant:
\begin{eqnarray*}
\emptyset=S^{j_0+1}\subseteq \tilde S^{j_0+1}\subseteq
S^{j_0}\subseteq \tilde S^{j_0}\subseteq\dots\\
\dots\subseteq S^3\subseteq
\tilde S^3\subseteq S^2 =\tilde S^2\subseteq S^1=\tilde S^1 =
X^{ss}
\end{eqnarray*}
\end{thm}

\begin{proo}
See \cite{C} (thm 2.1).
\end{proo}

\vskip .7 cm

\begin{remark} In general $S^i\neq\tilde S^i$: let us consider, for
  instance, $n=m=3$ and
$$A=\begin{pmatrix} 0 &0& x_0& x_1& x_2 \cr
    x_0& x_1& 0& 0& x_3\end{pmatrix}.$$
Since $D(A)=\{(0:0:t_1:t_2)\}\simeq \PP^1$, $A\in\tilde S^3$; but
$S^3=\emptyset$ (see also prop. \ref{codimensione}).
\end{remark}

\vskip 1 cm

\begin{corol}
If $m$ is odd and $A\in X^s=X^{ss}$ then $\codim D(A)\ge \frac {m+1}
2$.

If $m$ is even and $A\in X^{ss}$ (resp. $X^{s}$) then $\codim D(A)\ge
\frac m 2$ (resp. $>$).
\end{corol}

\begin{proo}
It suffices to notice that the previous theorem implies that $\tilde
S^{j_0+1}=\emptyset$ and that $S^{j_0}$ is the set of the properly
semi-stable points of $X$.
\end{proo}

\vskip .5 cm

\begin{prop}\label{codimensione}
If $m$ is odd, $A\in X$ is stable and $\codim D(A) = {\frac {m+1} 2}$,
  then, up to the action of
$\SL(I)\times\SL(W)\times\SL(V)$,
we have
$$A \simeq \begin{pmatrix}x_0 &\dots & x_{t-1} & 0 &\dots &0 &x_t\\ 0
  &\dots &0 &x_0 &\dots &x_{t-1} &x_{t+1}
\end{pmatrix},$$
where $t={\frac {m+1} 2}$.
\end{prop}
\begin{proof}
By the proof of theorem \ref{filtrazione} we have that for any
$\omega\in I$, $\dim(\omega\otimes V)\cap T = t$,
where $T$ is the image of $A$ as a subspace of $I\otimes V$.

Thus we have, up to a base change,
$$A \simeq \begin{pmatrix}x_0 &\dots & x_{t-1} & 0 &\dots &0 &x_t\\ 0
  &\dots &0 &y_0 &\dots &y_{t-1} &y_{t+1}
\end{pmatrix},$$
where $x_0,\dots, x_t$ and $y_0,\dots, y_t$ are linearly independent.

It is easily checked that $D(A) = V(x_0,\dots, x_t)\cup V(y_0, \dots,
y_t)\cup V(x_0,\dots,x_{t-1},
y_0,\dots,y_t)$ and since $\codim D(A) = t$, it must be $\codim V(x_0,\dots,x_{t-1},
y_0,\dots,y_t)= t$: this implies that $<x_0,\dots, x_{t-1}>=<y_0,\dots,y_{t-1}>$.

Moreover $x_t \neq a y_t$ for any $a\in \CC$ otherwise $A$ cannot be
stable.
\end{proof}

\begin{remark} The matrix above can exist if $n+1\ge t+1 = \frac
  {m+3} 2$, i.e. if $m\le 2 n -1$.

Since $A:W\hookrightarrow I\otimes V$ is injective, it must be $m+2\le
2 (n+1)$, i.e. $m\le 2 n$: thus in the odd case, the two requirements
are equivalent.
\end{remark}

\begin{corol}
Let $V_i=X^{ss}\setminus S^i$ e $\tilde V_i=X^{ss}\setminus \tilde
S^i$.

Then such subsets define a $G-$invariant increasing filtration:
$$\emptyset=V_1=\tilde V_1\subseteq V_2= \tilde V_2\subseteq
\tilde V_3\subseteq V_3\subseteq\dots$$
$$\dots\subseteq\tilde V_{j_0}\subseteq
V_{j_0}\subseteq \tilde V_{j_0+1}\subseteq V_{j_0+1}=X^{ss}.$$
In particular $V_2$ is the set of matrices that define vector bundles
and $V_{j_0}$ is the open set of the stable points in $X$.
\end{corol}

\vskip .6 cm

\begin{remark}
If $n$ is odd then  $V_{j(m)}=X^s=\tilde V_{j(m)+1} = V_{j(m)+1}=X^{ss}$.

Otherwise if $m$ is even then  $S^{j(m)}//G \simeq S^2\GG(\frac m 2,\PP^n)$
(theorem \ref{sottoinsieme}).
\end{remark}

\vskip .6 cm

All these results are needed to prove the following theorem:

\begin{thm}\label{stabilita}
Let $k=2$ and $m\in\NN$ odd.
$A\in \GG(m+2,I\otimes V)$ is G.I.T. stable if and only if $\mathcal F_A$ is $\mu$-stable.
\end{thm}

\begin{proo} See \cite{C} (thm. 3.1).
\end{proo}

Theorem \ref{main.thm} is a direct consequence of this equivalence
within the stability of the maps and the stability of the
cokernels.

\vskip 2 cm

%%%%%%%%%%%%%%%%%%%%%%%%%%%%%%%%%%%%%%%%%%%%%%%%%%%%%%%%%%%%
\section{Dimension of $S^j/G$}
%%%%%%%%%%%%%%%%%%%%%%%%%%%%%%%%%%%%%%%%%%%%%%%%%%%%%%%%%%%%

For any $j<j(m)$ we calculate the dimension of $S^j/G\subseteq
\mathcal M_{n,m,2}$ and we show that it is irreducible. In particular we
show that $S^2/G$
is the irreducible hypersurface that parameterizes all the sheaves in $
\mathcal M_{n,2}$ that are not bundles or, on the other hand, all the
points added to compactificate the moduli space of the rational normal
curves in $\PP^n$.

\

We remind that:
$$S^j=\{A\in X^{ss}|\exists ~0\neq \omega\in I \text{ such that }
\dim(T_A\cap R_\omega)\ge j+m-n\}.$$
Thus, if $j<j(m)$,
$$\frac{S^j}{\SL(W)} \simeq \{T\in \GG(m+1,\PP(I\otimes V))^{ss}|
\exists ~\omega\in I^* :\dim(T\cap
\PP(R_\omega))\ge j+m-n-1\}.$$
Let us define the incidence correspondence ${\mathcal I}_j\subseteq
\GG(m+1,\PP(I\otimes V))\times \PP(I)$ as:
$${\mathcal I}_j=\{(T,[\omega])|T\in  \GG(m+1,\PP(I\otimes V))^{ss},
[\omega]\in \PP(I), \dim(T\cap \PP(R_\omega))\ge j+m-n-1\} $$
and let $p_1$ and $p_2$ be the respective projections.
Since $S^1=X^{ss}$, we can suppose $2\le j<j(m) $.
Let us fix $[\omega]\in \PP(I)$: then
$$p_2^{-1}([\omega])\simeq\{T\in\GG(m+1,\PP(I\otimes V))^{ss}|\dim(T
\cap \PP(\omega\otimes V))\ge j+m-n-1\}$$ and:
\begin{equation*}
  \begin{split}
    \dim p_2^{-1}([\omega])&=
    (n+1-(j+m-n))(j+m-n)+\\
    &\qquad+(2(n+1)-(m+2))~(m+2-(j+m-n))=\\
    &=2mn-m^2 + 3n-m +(n-m)j + j - j^2.
  \end{split}
\end{equation*}
Hence ${\mathcal I}_j$ is irreducible (see \cite{Harris}, theorem 11.14)
of dimension $2mn-m^2 + 3n-m +1+(n-m)j + j - j^2$.

Now, if  $T\in p_1(\mathcal I_j)$ is a generic point, $p_1^{-1}(T)$ is
discrete, i.e. $\dim p_1^{-1}(T)=0$ that implies:
$$\dim S_{j}/\SL(W)=\dim p_1({\mathcal I}_j)=2mn-m^2 + 3n-m +1+(n-m)j
+ j - j^2.$$
Furthermore $S^j/\SL(W)$ is irreducible.

Since all the points of $S^j$ are stable under the action of $G$
(we are supposing $j<j(n)$), theorem \ref{main.thm} implies $$\dim
(S^j/G) = \dim S^j - \dim G = \dim (S^j/\SL(W)) - \dim \SL(I).$$

Hence we have:
\begin{thm}
$S^j/G$ is irreducible of codimension $(j+m-n)(j-1)-1$
for any $2\le j< j(m)$.
\end{thm}

In particular:
\begin{corol}If $n=m$ (boundary format)
$S^2/G$ is an irreducible hypersurface of $\mathcal M_{n,2}$
such that $$\mathcal M_{n,2} \setminus (S^2/G)\simeq S_{n,2}.$$
\end{corol}

By theorem \ref{sottoinsieme}, we know that, if $m$ is even,
the variety $\mathcal M_{n,m,2}\setminus (S^{j(m)}//G)$ is isomorphic to
$S^2\GG\left(\frac m 2, \PP(V)\right)$
and thus it is irreducible of dimension $(n-\frac m 2)(\frac m 2 +1)$,
i.e. the G.I.T. quotient $S^{j(m)}//G$ is of codimension
 $(n-\frac m 2)(\frac m 2 +1)$.

If $m$ is odd, then $S^{j(m)}=\emptyset.$

\vskip 2 cm

%%%%%%%%%%%%%%%%%%%%%%%%%%%%%%%%%%%%%%%%%%%%%%%%%%%%%%
\section{A torus action on $\mathcal M_{n,m,2}$}
%%%%%%%%%%%%%%%%%%%%%%%%%%%%%%%%%%%%%%%%%%%%%%%%%%%%%%

In the following two sections we compute the Euler characteristic of
$\mathcal M_{n,m,2}$ and an implicit formula for its Hodge numbers.
For this purpose, we will use the technique of Bialynichi-Birula
\cite{B}, that is based on the study of the action of a torus on a
smooth projective variety: such method was extensively used in the
last decade to compute the Betti numbers of smooth moduli spaces (see
for istance \cite{Kl}).

In fact let an algebraic torus $T$ act on a smooth projective variety
$Z$ and let $Z^T$ be its fixed points set. Then the Euler characteristics
of $Z$ and $Z^T$ are equal. Furthermore if $T=\CC^*$ is
1-dimensional, then all the cohomology
groups of $Z$ and their Hodge decomposition may be reconstructed from
the Hodge structure of the connected components $Z_i^T$ of
$Z^T$. In order to do that,  we fix a point $z_i\in Z_i^T$ for any
component and
we consider the action of $T$ on the tangent space $T_{z_i}Z$: let
$n_i$ be the number of positive weights of $T$ acting on $T_{z_i}Z$,
then we have:
\begin{thm}[Bialynichi-Birula]\label{BB}
There is a natural isomorphism:
$$\HH^{p,q}(Z) = \bigoplus_i \HH^{p-n_i,q-n_i}(Z^T_i).$$
\end{thm}
\begin{proof}
See \cite{B} and \cite{G}.
\end{proof}

\vskip .8 cm

Thus let us consider now the action of $T=\CC^*$ on $\mathcal
M_{n,m,2}$ defined by the morphism
$\rho:\CC^*\rightarrow  \GL(V)$ with weights
$c=(1,2,2^2,\dots,2^n)$: this choice is motiveted by the fact that
\begin{equation}\label{ci}
c_i-c_j=c_{i'}-c_{j'}\qquad \text{ if and only if } \qquad i = i'\text{ and } j= j'
\end{equation}
that will be useful later on.

For any $t\in \CC^*$, we will write $t(\cdot)$ to denote the image of
$\cdot$ by the map $\rho(t)$.

Let $A=\begin{pmatrix} f_0 \dots f_{m+1} \\ g_0 \dots
g_{m+1}\end{pmatrix}\in \mathcal M_{n,m,2}$ be a fixed point then
$$t(A)=\begin{pmatrix} t(f_0) \dots t(f_{m+1}) \\ t(g_0) \dots
t(g_{m+1})\end{pmatrix}\sim A$$
for any $t\in \CC^*$. Thus it is defined a morphism $\tilde
\rho:\CC^*\rightarrow \Aut(I)\times \Aut(W)$, such that
$\rho(t)(A)=\tilde\rho(t)(A)$ for any $t\in \CC^*$.

Thus for any fixed point $A$,  $\rho$ induces an action of $\CC^*$
on $I$ and $W$: let $P(t)$ and $Q(t)$ be the components of $\tilde
\rho$ in $\Aut(I)$ and $\Aut(W)$ respectively, then
$t(A)=P(t) ~ A ~ Q(t)^{-1}$ for any $t$ in $\CC^*$.
We can suppose that such action is diagonal and that
it is defined by the weights $(a_0,a_1)$ and $(b_0,\dots, b_{m+1})$
respectively (at the moment we do not fix any order for such weights,
we will do it later on).

If $f_k = \sum r_i x_i$ then $\sum r_i t^{c_i} x_i=t(f_k)=\sum r_i
t^{a_0 - b_k}x_i$, and
since $c_i\neq c_j$ if $i\neq
j$, it must be $f_k =
r_{i_k}x_{i_k}$ for a suitable $i_k\in\{0,\dots,n\}$ and with
$r_{i_k}\in \CC$; moreover it results $a_0 - b_k = c_{i_k}$
for any   $k$ such that $r_{i_k}\neq 0$.

Similarly we have $g_k=s_{j_k}x_{j_k}$ with $s_{j_k}\in \CC$,
$j_k\in\{0,\dots,n\}$ and $a_1 - b_k = c_{j_k}$ for any $k$ such that
$r_{j_k}\neq 0$.

Thus the matrix $A$ is monomial with respect to the bases of $I$
and $W$ chosen. Moreover the weights $(a_0,a_1)$ and
$(b_0,\dots,b_{m+1})$ are the solution of a system:
\begin{equation}\label{spf}
\begin{cases}
a_0 - b_k = c_{i_k} \quad \forall ~k \text{ s.t. } f_{k}\neq 0\\
a_1 -  b_k = c_{j_k} \quad \forall ~k \text{ s.t. } g_{k}\neq 0\\
\end{cases}
\end{equation}

Since $A$ is stable, there exists $\tilde k$ such that $f_{\tilde
k},g_{\tilde k}\neq 0$, thus, by (\ref{spf}), it follows  that $a_0-a_1 =
c_{i_{\tilde k}} -c_{j_{\tilde k}}$: it is easy to check that  if
(\ref{spf}) admits  a solution, then such solution is unique up to an
additive constant; for this reason we can suppose $a_0=0$.

Now we can fix an order on the base of $W$ chosen (we did
not do it before): in fact we can suppose $f_k=0$ if and only if  $k>
k_0$ where $k_0\in\{1,\dots,m+1\}$: moreover we can take $b_0\ge b_1
\ge\dots\ge b_{k_0}$ and, if $k_0\le m$, we can also take
$b_{k_0+1}\ge b_{k_0+2}\ge\dots \ge b_{m+1}$.
In particular we have $c_{i_0}\le c_{i_1}\le \dots \le c_{i_{k_0}}$ and
$c_{j_{k_0+1}}\le\dots\le c_{j_{m+1}}$, that implies $i_0\le i_1\le
\dots\le i_{k_0}$ and $j_{k_0+1}\le\dots\le j_{m+1}$.

Let $k_1,\dots,k_z\le k_0$ be such that $f_{k_j},g_{k_j}\neq 0$ for
any $j=1,\dots, z$: it must be $z\ge 1$ and
$a_1=c_{j_{k_1}}-c_{i_{k_1}}$.
Thus (\ref{spf}) becomes:
\begin{equation}\label{spf2}
\begin{cases}
b_k= -c_{i_k} \quad &\forall ~k\le k_0\\
b_k=a_1 - c_{j_k} \quad &\forall ~k> k_0\\
a_1=c_{j_{k_s}}-c_{i_{k_s}}\quad &\forall ~s=1,\dots,z
\end{cases}
\end{equation}
By (\ref{ci}) and since
$(f_{k_j},g_{k_j})\neq (f_{k_1},g_{k_1})$ if $s=2,\dots, z$, we can
either suppose $z=1$ or $a_1=0$ that implies $i_{k_s}=j_{k_s}$ for any
$s=1,\dots, z$.
Thus we have to distinguish two cases:
\begin{enumerate}
\item $a_1\neq 0$, $z=1$
\item $a_1=0$, $z\ge 1$
\end{enumerate}

Under each of these hypothesis, it is easy to show that the system
(\ref{spf2}) admits a unique solution that defines a fixed point
$A\in\mathcal M_{n,m,2}$ by the action of
$\rho$.

In order to have a total description of the fixed points, we will
consider each case separately:

\begin{enumerate}

\item Let us define
\begin{equation}\label{primo.caso}
A_{I,J}=\begin{pmatrix}
 x_{i_0} & x_{i_1} & \dots & x_{i_t} & 0       & \dots & 0 \\
 x_{j_0} & 0       & \dots & 0       & x_{j_1} & \dots & x_{j_t}
\end{pmatrix}
\end{equation}
where $I=(i_0,\dots,i_t)$ and $J=(j_0,\dots,j_t)$, with
$i_1<\dots<i_t$, $j_1<\dots<j_t$, $i_0<j_0$ and $i_0\neq i_s$,
$j_0\neq j_s$ for any $s=1,\dots,t$.
%%%%%%%%%%%%%%%%%%%%%%%%%%%%%%%%%%%%%%%%%%%%%%%%%%%?? controllare

\noindent It is easy to see that under this assumption the matrices
$A_{I,J}$'s  are stable and  determine uniquely all the fixed point of
$\rho$ with $a_1\neq 0$.

\vskip 5 mm

\item The matrices fixed by $\rho$ with $a_1=0$ are given by
\begin{equation}\label{secondo.caso}
A_{\omega}^i = (\omega_1~x_{i_1},\dots,\omega_{m+2} ~x_{i_{m+2}})
\end{equation}
with $\omega=(\omega_1,\dots,\omega_{m+2})\in I^{m+2}$ and
$i=(i_1,\dots,i_{m+2})$ where  $0\le i_1\le\dots\le i_{m+2}\le n$.

\noindent Since $\dim I=2$, if $i_j=i_{j+1}$ then we can suppose
that $\{\omega_{i_j}, \omega_{i_{j+1}}\}$ is the base of $I$ fixed above.
Moreover there cannot exist a $j$ such that $i_{j-1}=i_j=i_{j+1}$
otherwise $A_{\omega}^i$ cannot be stable.
Thus, in particular,
\begin{equation}\label{li}
l(i)=\#\{i_j|i_j\neq i_k \text{ for any } k\neq j\}
\end{equation}
is odd.

\noindent It is easy to check that $A_\omega^i$ and $A_{\omega'}^{i'}$
are contained in the same connected component of $\mathcal M_{n,m,2}^\rho$
if and only if $i=i'$ and that such component is isomorphic to
$\PP(I)^{l(i)}/\SL(2)$. In particular $A_\omega^i$ is stable if and
only if  the corresponding poing in $\PP(I)^{l(i)}$ is stable under
the action of $\SL(2)$. The stable points of $\PP(I)^{l(i)}$
are described by:
\end{enumerate}

\vskip .6 cm

\begin{prop}Let $l\in \NN$ be odd and let the group
$\SL(I)$ act on $Y=\PP(I)^l=\PP(I)\times\dots\times\PP(I)$;
for any $\omega\in Y$, let $I_k(\omega) =
\{j\in(1,\dots,n)|\omega_i=\omega_k\}$;
then
$$Y^s=Y^{ss}=\{\omega \in Y|~ \# I_k(\omega) < \frac l 2 \quad \text{
  for any } k=1,\dots,l\}.$$
\end{prop}

\begin{proo}
It is a direct consequence of the Hilbert-Mumford criterion for stability.
(see also \cite{MFK}).
\end{proo}

\vskip .8 cm

The Hodge numbers of $M_l=\PP(I)^l/SL(I)$ are given by the following:
\begin{thm}\label{hodge.ml}
Let $l$ be odd. Then:
$$\hh^{p,q}(M_l)=\begin{cases} 0\qquad &\text{if } p\neq q\\
      1 + (l-1) + \dots + \begin{pmatrix} l-1\\
\min(p,l-3-p)\end{pmatrix}\qquad &\text{if } p=q \end{cases}.$$
In particular, the Poincar\'e polynomial is:
$$P_t(M_l)=1 + \hh^{1,1}t^2 + \dots +  \hh^{j,j}t^{2j}+\dots+t^{2l-6}.$$
and the Euler characteristic is given by:
$$\e(M_l)= \sum_{p=0,\dots, l-3} \hh^{p,p}$$
\end{thm}

\begin{proo}
See \cite{K} pag. 193.
\end{proo}

\vskip 5 mm

By the classification of the fixed points of $\mathcal M_{n,m,2}$, we
thus have:
\begin{corol}
$h^{p,q}(\mathcal M_{n,m,2})=0$ for any $p\neq q$.
\end{corol}

\vskip 5 mm

We are now ready  to compute the Euler characteristic of $\mathcal M_{n,m,2}$:

\begin{thm}\label{eul.char}
Let $m$ be odd and let $t=\frac{m+1} 2$. Then the Euler characteristic
of  $\mathcal M_{n,m,2}$ is given by:
\begin{equation}\label{eul}
\e(\mathcal M_{n,m,2}) = \begin{pmatrix} n+1\\ 2  \end{pmatrix}
           \begin{pmatrix} n \\t\end{pmatrix}^2
+\sum_{d=1}^{n-t}  \begin{pmatrix}n+1\\ t-d\end{pmatrix}
\begin{pmatrix} n+1 -t +d\\ 2d +1 \end{pmatrix}
  \e(\PP(I)^{2d+1}/\SL(I))\enspace .
\end{equation}
\end{thm}
\begin{proo}
By theorem \ref{BB}, it results $\e(\mathcal M_{n,m,2})=\sum_i \e
(\mathcal M_{n.m,2}^T)_i$, where $(\mathcal M_{n.m,2}^T)_i$ are the
connected components of the fixed point of $\mathcal M_{n,m,2}$ under
the action of the torus $T$ considered.

The points $A_{I,J}$, defined in (\ref{primo.caso}), represent
discrete components of such space and since they are uniquely determined by
$I=(i_0,\dots,i_t)$ and $J=(j_0,\dots,j_t)$, with
$i_1<\dots<i_t$, $j_1<\dots<j_t$ and $i_0<j_0$,
it is easy to compute that they are exactly
\begin{equation}\label{eul1}
\begin{pmatrix} n+1\\ 2  \end{pmatrix}
\begin{pmatrix} n \\t\end{pmatrix}^2.
\end{equation}

On the other hand, the matrices $A^i_\omega$, defined in
(\ref{secondo.caso}), form connected components
determined by $i=(i_1,\dots,i_{m+2})$ and isomorphic to $M_{l(i)}$, where
$l(i)$ is defined in (\ref{li}).

Let $d(i)=m+2-l(i)$ the number of the couples of equal terms in $i$:
for any $d\ge 1$ the number of the admissible vectors
$i=(i_1,\dots,i_{m+2})$ with $d(i)=d$ is given by
\begin{equation*}
\begin{pmatrix}n+1\\d\end{pmatrix} ~ \begin{pmatrix}n+1-d\\m+2-2d\end{pmatrix}
\end{equation*}
thus the Euler characteristic of the set of such matrices is given by
\begin{equation}\label{eul2}
\begin{split}
\sum_{d=m-n+1}^{t-1} \begin{pmatrix}n+1\\d\end{pmatrix} ~
 \begin{pmatrix}n+1-d\\m+2-2d\end{pmatrix} ~
  \e(M_{m+2-2d}) \\
=\sum_{d=1}^{n-t}\begin{pmatrix}n+1\\ t-d\end{pmatrix}
\begin{pmatrix} n+1 -t +d\\ 2d +1 \end{pmatrix}
   ~ \e(M_{2d+1}).
\end{split}
\end{equation}
(\ref{eul}) is obtained by summing (\ref{eul1}) with (\ref{eul2}).
\end{proo}
\vskip 2 cm

%%%%%%%%%%%%%%%%%%%%%%%%%%%%%%%%%%%%%%%%%%%%%%%%%%%%%%
\section{Betti numbers}
%%%%%%%%%%%%%%%%%%%%%%%%%%%%%%%%%%%%%%%%%%%%%%%%%%%%%%

We compute the numbers $n_i$ for any fixed point in a connected
component $(\mathcal M_{n,m,2}^T)_i$. We remind that $n_i$ represents
the number of positive weights of $T=\CC^*$ acting on the tangent
space of $\mathcal M_{n,m,2}$ at the fixed points.
These numbers will yield to the computation of  the Betti numbers
of  $\mathcal M_{n,m,2}$ for any odd $m$.

In particular we get a topological description of the moduli space of
the rational normal curves on $\PP^n$ for any odd $n$.

\vskip 6 mm

Let $A\in \mathcal M_{n,m,2}$ be a fixed point for $\rho$.
Then $\rho$ induces an action on the tangent space $T_A \mathcal
M_{n,m,2}$. By theorem \ref{main.thm}, such vector space is isomorphic
to the tangent space of the Maruyama scheme
$\mathcal M_{\PP^n}(m;c_1,\dots,c_n)$ at the point corresponding to
the sheaf $\mathcal F_A$ and thus it is isomorphic to $\Ext^1(\mathcal
F_A,\mathcal F_A)$ (see \cite{Mar1} and \cite{Mar2}).

By the sequence (\ref{suc.esatta}) that defines the sheaf $\mathcal
F_A$, it is easily checked that  $\Ext^1(\mathcal F_A,
\mathcal F_A)$ is contained in the exact sequence:
$$0\to \Hom(\mathcal F_A,\mathcal F_A)\to \Hom(W\otimes
\OPN(1),\mathcal F_A) \to \Hom(I\otimes \OPN(1),\mathcal F_A) \to
\Ext^1(\mathcal F_A,\mathcal F_A)\to 0.$$
Moreover $\HH^0(\mathcal F_A)=(W\otimes V)/a(I)$ where $a:I\hookrightarrow
W\otimes V$ is the map induced by $A^*:I\otimes \OPN\hookrightarrow
W\otimes \OPN(1)$ and $\HH^0(\mathcal F(-1))=W$.

Thus, it results $\Hom(W\otimes\OPN(1),\mathcal F_A)=W^*\otimes
\HH^0(\mathcal F_A(-1)) = W^*\otimes W$ and $\Hom(I\otimes \OPN, \mathcal
F_A)=I^*\otimes (W\otimes V)/a(I)$.

In particular the weights of the action $\rho$ on $\Ext^1(\mathcal
F_A,\mathcal F_A)$ are easily computed using the sequence:
\begin{equation}\label{pesi.ext}
0\to \CC \to W^*\otimes W \to I^*\otimes \frac {W\otimes V} {a(I)} \to
\Ext^1(\mathcal F_A,\mathcal F_A)\to 0.
\end{equation}

In the previous section we have seen that for any fixed matrix $A\in X$, $\rho$
induces an action on $I$ and $W$ defined by the weights $(a_0,a_1)$ and
$(b_0,\dots, b_{m+1})$ described by (\ref{spf2}), where, we remind,
$c=(1,2,\dots,2^n)$.

For any $A\in (\mathcal M_{n,m,2}^T)_i$ we write $n(A)$ in place of
$n_i$ and moreover we define $n_1(A)$ as the number of the positive
weights of $\rho$ on $W^*\otimes W$ and similarly $n_2(A)$ as the number
of the positive weights on $I^*\otimes (W\otimes V)/a(I)$.
Thus by the sequence (\ref{pesi.ext}), it results $n(A)=n_2(A)-n_1(A)$.

\vskip 5 mm

In order to calculate $n_1(A)$ and $n_2(A)$ for all the fixed matrices
by the action of $\rho$, we need to distinguish the cases
described above:

\begin{prop}\label{ni}
\
\begin{enumerate}
\item Let $A_{I,J}$ be defined as in (\ref{primo.caso}); then:
\begin{align*}
n_1(A_{I,J}) = ~ & 4tn +2t+2n -1 - \sum_{s=0}^t i_s-\sum_{s=0}^t j_s -
\sum_{i_s>i_0} i_s - \sum_{j_s>i_0,s\ge 1} j_s \\
& - \#\{s=1,\dots,t|j_s>j_0\} -
i_0\cdot\#\{s=1,\dots,t|j_s\le i_0\} \\
\intertext{and}\\
n_2(A_{I,J}) = ~ & \begin{pmatrix} m+2 \\ 2 \end{pmatrix}.\\
\\
\intertext{\item Let $A_\omega^i$ be defined as in
(\ref{secondo.caso}); then:}\\
n_1(A_\omega^i) = ~ &2(m+2)n -2\sum_{s=0}^{m+1}i_s~\\
\intertext{and}\\
n_2(A_\omega^i) = ~& \begin{pmatrix}m+2 \\ 2\end{pmatrix} +
\frac{m+2-l(i)} 2.
\end{align*}

\end{enumerate}
\end{prop}
\begin{proof}
It is just a direct computation.
\end{proof}

\vskip 1 cm

Proposition \ref{ni} and theorem \ref{hodge.ml} give us the right
ingredients to apply theorem
\ref{BB} of Bialynichi-Birula. Thus we have an algorithm to compute
the Betti numbers of $\mathcal M_{n,m,2}$ for any $m\ge n$, and in
particular of $\mathcal M_{n,n,2}$ the compactification of the variety
$S_n$ of the rational normal curves.

In fact, let $b_i(n)=\dim\HH^{i}(\mathcal M_{n,n,2},\QQ)$: the
following table provides the values of $b_i(n)$, for $n=2,3,5,7$ and
for all the even $i=0,\dots, 36$.

\vskip 1 cm
\begin{tabular}{|c||c|c|c|c|c|c|c|c|c|c|c|c|c|c|c|c|c|c|c|c} \hline
 $n$ &
 $ b_0    $ & $ b_2    $ & $ b_4    $ & $ b_6    $ & $ b_8    $  &
 $ b_{10} $ & $ b_{12} $ & $ b_{14} $ & $ b_{16} $ & $ b_{18} $  &
 $ b_{20} $ & $ b_{22} $ & $ b_{24} $ & $ b_{26} $ & $ b_{28} $  &
 $ b_{30} $ & $ b_{32} $ & $ b_{34} $ & $ b_{36} $
 \\\hline\hline
$ 2$ & 1 & 1 & 1 & 1 & 1 & 1 &  &  &  &  &  &  &  & & & & & &
\\\hline
$ 3$ & 1 & 1 & 3 & 4 & 7 & 8 & 10 & 8  & 7  & 4 &3 & 1 & 1 & & & & & &
\\\hline
$ 5$ & 1 & 1 & 3 & 4 & 8 &11 & 18 & 24  & 35  & 45 & 61 & 74 & 93
&106 & 122 & 128 & 134 & 128 & 122
\\ \hline
$ 7$ & 1 & 1 & 3 & 4 & 8 &11 & 19 & 26  & 40  & 54 &77 & 100 & 134 &165 &205 &
 242 & 289 & 334 & 400
\\ \hline
$ 9$ & 1 & 1 & 3 & 4 & 8 &11 & 19 & 26  & 41  & 56 &82 &110 & 154 &202 &273 &
352 & 461 & 595 & 750
\\ \hline
\end{tabular}

\vskip 1 cm

See also \cite{EPS}, for the computation of the Betti numbers of
$\mathcal M_{3,3,2}$.

\vskip 5 mm

By this table, it seems that, for any $i\ge 0$ and $n>>0$,
the value of $b_i(n)$ is costant. In particular we have:

\begin{prop}
$b_2(n)=h^{1,1}(\mathcal M_{n,n,2})=1$ for any odd $n$.
\end{prop}

\begin{proof} By prop. \ref{ni}, it follows that $n_i(A_\omega^i)\ge
  3$ for any $A_\omega^i$ defined as in (\ref{secondo.caso}).
Thus, by theorem \ref{BB}, the points represented by the matrices
$A_\omega^i$, do not give any contribute to $b_2(n)$.

Moreover it is not difficult to see that, if $n>3$, the only matrix
$A_{I,J}$, as in (\ref{primo.caso}), such that $n_i(A_{I,J})=1$ is given by
$I=(n-t,n-t+1,\dots,n)$ and $J=(n-t-2,n-t+1,n-t+2,\dots,n)$, where $t=
\frac{n+1} 2$.
\end{proof}

\vskip 5 mm
\begin{remark}
It would be interesting to have a description of the Chow rings of
$\mathcal M_{n,m,2}$: in \cite{Sc}, the author studies the Chow ring of
the Hilbert compactification $\mathcal H_3$ of the moduli space of the
twisted cubics in $\PP^3$.
\end{remark}

\vskip 2 cm

%%%%%%%%%%%%%%%%%%%%%%%%%%%%%%%%%%%%%%%%%%%%%%%%%%%%%%
\providecommand{\bysame}{\leavevmode\hbox to3em{\hrulefill}}

\end{document}